\documentclass[11pt,english,reqno]{amsart}

\usepackage[paperwidth=210mm,paperheight=297mm,inner=3.1cm,outer=3.1cm,top=2.5cm,bottom=2.5cm]{geometry}

\usepackage[T1]{fontenc}
\usepackage{amssymb}

\usepackage{paralist}
\usepackage{babel}
\usepackage[pdftex]{graphicx}    
\usepackage{amsmath}
\usepackage{amsthm}

\newcommand\mbb{\mathbb}

\newcommand\mcal{\mathcal}

\newcommand\ol{\overline}

\newcommand\wt{\widetilde}

\newcommand\sC{\mcal{C}}

\newcommand\sR{\mcal{R}}

\newcommand\sT{\mcal{T}}

\newcommand\C{\mbb{C}}

\newcommand\R{\mbb{R}}
\newcommand\Z{\mbb{Z}}

\DeclareMathOperator*\Hom{Hom}

\DeclareMathOperator*\QM{QM}

\newcommand\isom{\cong}

\renewcommand\epsilon{\varepsilon}
\renewcommand\ge{\geqslant}

\renewcommand\le{\leqslant}

\renewcommand\phi{\varphi}

\renewcommand\theta{\vartheta}

\numberwithin{equation}{section}

\theoremstyle{plain}
\newtheorem{Thm}[equation]{Theorem}
\newtheorem{Prop}[equation]{Proposition}
\newtheorem{Cor}[equation]{Corollary}
\newtheorem{Lemma}[equation]{Lemma}

\newtheorem*{Thm*}{Theorem}
\newtheorem*{Prop*}{Proposition}
\newtheorem*{Cor*}{Corollary}
\newtheorem*{Lemma*}{Lemma}
\newtheorem*{Sublemma*}{Sublemma}
\newtheorem*{Conjecture*}{Conjecture}

\theoremstyle{definition}
\newtheorem{Def}[equation]{Definition}

\newtheorem{Examples}[equation]{Examples}

\newtheorem{Remark}[equation]{Remark}

\newtheorem*{Def*}{Definition}
\newtheorem*{Defs*}{Definitions}
\newtheorem*{Example*}{Example}
\newtheorem*{Examples*}{Examples}
\newtheorem*{Exercise*}{Exercise}
\newtheorem*{LemmaDef*}{Lemma and Definition}
\newtheorem*{Notation*}{Notation}
\newtheorem*{Problem*}{Problem}
\newtheorem*{Question*}{Question}
\newtheorem*{Remark*}{Remark}
\newtheorem*{Remarks*}{Remarks}
\newtheorem*{Warning*}{Warning}

\newtheorem*{Text*}{}

\newcommand\Co{\sC^0}
\DeclareMathOperator*\re{Re}

\begin{document}
\title{Positivity of continuous piecewise polynomials} \author{Daniel
  Plaumann} \address{Department of Mathematics, University of
  California, Berkeley, CA 94720, USA}
\email{plaumann@math.berkeley.edu} 

\date{\today}

\thanks{The author gratefully
  acknowledges support through a Feodor Lynen fellowship from the Alexander
  von Humboldt foundation.}
\maketitle

\begin{abstract}
  Real algebraic geometry provides certificates for the positivity of
  polynomials on semi-algebraic sets by expressing them as a suitable
  combination of sums of squares and the defining inequalitites. We
  show how Putinar's theorem for strictly positive polynomials on
  compact sets can be applied in the case of strictly positive
  piecewise polynomials on a simplicial complex. In the
  $1$-dimensional case, we improve this result to cover all
  non-negative piecewise polynomials and give explicit degree bounds.
\end{abstract}

\section*{Introduction}

Let $\Delta=\sigma_1\cup\cdots\cup\sigma_k$ be a simplicial complex in
$\R^n$ with vertices $v_1,\dots,v_m$. Let $\Co(\Delta)$ denote the
algebra of all continuous piecewise polynomials on $\Delta$,
consisting of all continuous functions $F\colon\Delta\to\R$ such that
the restriction of $F$ to each $\sigma_i$ is given by a polynomial. It
has been studied in connection with splines, where the functions are
also required to be differentiable to some order. For a good survey,
see Strang \cite{MR0327060} and references given there. A detailed
analysis of $\Co(\Delta)$ from the point of view of combinatorics and
commutative algebra is due to Billera \cite{MR1013666}.

The algebra $\Co(\Delta)$ has a beautiful description by generators
and relations, given in terms of its \emph{tent functions} or
\emph{Courant functions}: These are the unique piecewise linear
functions $T_i\colon\Delta\to\R$ such that $T_i(v_i)=1$ and
$T_i(v_j)=0$ if $i\neq j$. The tent functions generate $\Co(\Delta)$
and satisfy certain obvious relations, which Billera has shown to be
sufficient to completely describe $\Co(\Delta)$ (see
Thm.~\ref{thm:Billera} below). The relations are in fact identical to those of
the \emph{Stanley-Reisner ring} (or face ring) of $\Delta$, plus one
additional relation that accounts for the fact that the tent functions
sum to $1$ (see \cite{MR1013666}, \cite[Def.~1.6]{MR2110098}).

\smallskip In this paper, we show how the tent functions can also be
used to characterize positive and non-negative functions in
$\Co(\Delta)$. In general, much work in real algebraic geometry has
been concerned with so-called certificates for positivity: Let
$h_1,\dots,h_r\in\R[t]$ be real polynomials in $n$ variables
$t=(t_1,\dots,t_n)$, and let $S$ be the basic closed semi-algebraic
set $\{x\in\R^n\ |\ h_1(x)\ge 0,\dots,h_r(x)\ge 0\}$. The convex cone
$M\subset\R[t]$ generated by all polynomials $g^2$ and $g^2h_i$ is
called a \emph{quadratic module}. Putinar \cite{MR1254128} has shown
that $M$ contains all polynomials that are strictly positive on $S$ if
there exists $N\in\Z_+$ such that $N-\sum t_i^2\in M$, in which case
$M$ is called \emph{archimedean}. If $M$ is archimedean, then $S$ is
clearly compact. Schm{\"u}dgen's positivstellensatz \cite{MR1092173}
(which predates that of Putinar) can be rephrased as saying that the
quadratic module generated by all square-free products of
$h_1,\dots,h_r$ (called the preordering) is archimedean whenever $S$
is compact. These and many related results have attracted attention in
optimization because membership of a polynomial in a preordering or
quadratic module can be checked (in practice rather efficiently) by a
semidefinite program. Good general references on the subject are the
books of Marshall \cite{MR2383959} and Prestel and Delzell
\cite{MR1829790}.

\smallskip A simplicial complex $\Delta$ is, of course, a special kind
of semi-algebraic set, and the algebra $\Co(\Delta)$ can be
interpreted as the ring of polynomial functions on an affine algebraic
variety containing $\Delta$ (which is just the Zariski-closure of
$\Delta$ if the ambient dimension $n$ is large and the vertices of
$\Delta$ are in general position; see Prop.~\ref{prop:DirectLimit}). It is
therefore possible to translate known results into this setup, in
particular Putinar's positivstellensatz (Thm.~\ref{thm:PutPPF}):

\begin{Thm*}
  Let $\Delta\subset\R^n$ be a simplicial complex with $m$
  vertices. Let $T_1,\dots,T_m\in\Co(\Delta)$ be the tent functions on
  $\Delta$. If a function $F\in\Co(\Delta)$ is strictly positive on
  $\Delta$, then there exist sums of squares
  $S_i$ in $\Co(\Delta)$, $i=0,\dots,m$, such that
\[
F= S_0+\sum_{i=1}^m S_i T_i.
\]
\end{Thm*}

\noindent Beyond pointing out this application of a known result, the
contribution of this paper is a strengthening in the case when
$\Delta$ is $1$-dimensional:

\begin{Thm*}
  Let $\Delta\subset\R^n$ be a simplicial complex of dimension $1$
  with $e$ edges and $m$ vertices, of which $m_0$ are isolated. Let
  $T_1,\dots,T_m\in\Co(\Delta)$ be the tent functions on $\Delta$. A
  function $F\in\Co(\Delta)$ is non-negative on $\Delta$ if and only
  if there exist sums of squares $S, S_{ij}$ in $\Co(\Delta)$ such
  that
\[
F=S + \sum_{(i,j)\in E} S_{ij} T_iT_j.
\]
More precisely, there exist such $S$ and $S_{ij}$ of degree at
most $\deg(F)+6(e-1)+1$, where $S$ is a sum of at most $2e+m_0$
squares and each $S_{ij}$ is a sum of two squares.
\end{Thm*}

Positive polynomials on general semialgebraic subsets of real
algebraic curves have been studied extensively by Scheiderer in
\cite{MR2020709}. The criteria developed there cover all semialgebraic
subsets of irreducible curves. The proof of the main result above
given here is elementary and (in principle)
constructive. Alternatively, it is possible to obtain a more abstract
proof (without the degree bounds) using the local results and
local-global principle in \cite{MR2020709}. The general case of
semialgebraic subsets of a reducible curve is not settled
completely. Most of the results obtained by the author in
\cite{MR2652535} only apply to sums of squares. The existence of
degree bounds for quadratic modules (\emph{stability}) is also open in
most $1$-dimensional cases, even for irreducible curves (except when
the curve is rational or elliptic). The case of a simplicial complex
is thus a very particular one, not only because the curves involved
are just lines, but also because the generators $T_iT_j$ of the
quadratic module vanish in the intersection points (the vertices),
which turns out to be very helpful.

\smallskip \emph{Acknowledgements.} I would like to thank Mihai
Putinar, Claus Scheiderer, and Bernd Sturmfels for helpful discussions
on the subject of this paper. I am also grateful to the referee for
suggesting several improvements, including the use of
quadratic modules instead of preorderings and slightly better degree bounds.

\section{Preliminaries}\label{sec:prelim}

Let $v_0,\dots,v_d\in\R^n$ be affinely independent vectors, $d\le
n$. We will write $\sigma(v_0,\dots,v_d)$ for the $d$-simplex spanned
by $v_0,\dots,v_d$, which is the convex hull of
$\{v_0,\dots,v_d\}$. It is equivalent to the $d$-dimensional standard simplex
$\sigma(e_1,\dots,e_{d+1})$ up to an affine change of coordinates. The
faces of $\sigma(v_0,\dots,v_d)$ are precisely the $2^{d+1}$
subsimplices $\sigma(v_{i_0},\dots,v_{i_e})$ with
$\{i_0,\dots,i_e\}\subset\{0,\dots,d\}$, $e\le d$.

A \emph{simplicial complex} is a union
$\Delta=\sigma_1\cup\cdots\cup\sigma_k$, where each $\sigma_i$ is a
simplex in $\R^n$, such that $\sigma_i\cap\sigma_j$ is a face of both
$\sigma_i$ and $\sigma_j$, for all $i,j=1,\dots,k$. One can always
assume that $\sigma_i\nsubseteq\sigma_j$ for all $i,j$. In that case,
$\Delta$ is called \emph{pure dimesional} if all $\sigma_i$ have the
same dimension.

Let $d_i=\dim(\sigma_i)$, and write 
$\sigma_i=\sigma(v_{i0},\dots,v_{id_i})$. We denote by
$V_\Delta=\{v_{ij}\ |\ i=1,\dots,k, j=0,\dots,d_i\}$ the set of all
vertices of $\Delta$. The simplicial complex $\Delta$ is uniquely determined
by the set $V_\Delta$ together with the information which of the
simplices $\sigma(V')$ for $V'\subset V_\Delta$ are contained in $\Delta$.

A function $f\colon\Delta\to\R$ will be called a
\emph{piecewise polynomial on $\Delta$ }if there exist polynomials
$g_1,\dots,g_k\in\R[t_1,\dots,t_n]$ such that
$f|_{\sigma_i}=g_i|_{\sigma_i}$, for all $i=1,\dots,k$. We will study
positivity in the algebra
\[
\Co(\Delta)=\bigl\{f\colon\Delta\to\R\ \bigl|\ f\text{ is continuous
  and piecewise polynomial}\bigr\}.
\]

As explained in the introduction, we will work with two different
descriptions of $\Co(\Delta)$, the first of which is given as the
coordinate ring of an affine $\R$-variety. For the purpose of this
paper, such a variety is given by an $\R$-algebra $\R[U]$ that is
reduced (i.e.~without nilpotents) but not necessarily integral. The
geometric object $U$ associated with $\R[U]$ is the set
$\Hom(\R[U],\R)$ of $\R$-algebra homomorphisms. After fixing
coordinates, i.e.~a surjection $\phi\colon \R[t_1,\dots,t_N]\to\R[U]$
for some $N$ and a finite set of generators $f_1,\dots,f_r$ of
$\ker(\phi)$, the set $U$ is in canonical bijection with $\{x\in\R^N\
|\ f_1(x)=\cdots f_r(x)=0\}$, the common zero-set of
$f_1,\dots,f_r$. Note that we only consider real points here (i.e.~we
identify $U$ with $U(\R)$ in the notation of \cite{MR2020709},
\cite{MR2652535}), which is sufficient for our needs.

\begin{Prop}\label{prop:DirectLimit}
  Let $\Delta=\sigma_1\cup\cdots\cup\sigma_k\subset\R^n$ be a simplicial complex. For all $i\le j=1,\dots,k$, let $U_{ij}$ be the affine hull of $\sigma_i\cap\sigma_j$. Let $U$ be the affine $\R$-variety obtained as the direct limit of the directed system $\{U_{ij}\}$ ordered by the inclusions $U_{ij}\to U_{ii}$, $U_{ij}\to U_{jj}$. Then 
\[
\Co(\Delta)\isom\R[U]\isom\biggl\{(f_1,\dots,f_k)\in\prod_{i=1}^k\R[U_{ii}]\ \bigl|\ f_i|_{U_{ij}}=f_j|_{U_{ij}}\text{ for all }i\neq j\biggr\}.
\]
\end{Prop}

\begin{proof}
  It suffices to note that the values of a
polynomial $g_i$ on $\sigma_i\cap\sigma_j$ uniquely determine its values on the
Zariski-closure, which is $U_{ij}$. 
\end{proof}

If all vertices of $\Delta$ are in sufficiently general position, so that the affine spans $U_{ii}$ of $\sigma_i$ are all distinct subspaces of $\R^n$, then $U$ is isomorphic to the union $\bigcup_{i=1}^k U_{ii}$. One can always arrive at this situation by embedding $\Delta$ into a higher-dimensional ambient space, if necessary.

Once we fix isomorphisms $U_{ii}\isom\R^{d_i}$ (even though there is no
canonical choice), this fixes $U_{ij}$ as affine subspaces of $U_{ii}$, $U_{jj}$, and we can think of the elements of $\Co(\Delta)$ as
$k$-tuples of polynomials together with a compatibility condition:
\[
\Co(\Delta)\isom\bigl\{(g_1,\dots,g_k)\in\prod_{i=1}^k R[t_1,\dots,t_{d_i}]\
\bigl|\ g_i|_{U_{ij}}=g_j|_{U_{ij}}\text{
  for all }i,j=1,\dots,k\bigr\}.
\]

The second, more intrinsic description of $\Co(\Delta)$ in terms
of its tent functions, was given by Billera in \cite{MR1013666}: Let
$V_\Delta=\{v_1,\dots,v_m\}$ be the set of vertices in $\Delta$, as
above. For $i=1,\dots,m$, let $T_i$ be the \emph{tent function} or
\emph{Courant function} on $\Delta$, which is the unique piecewise
linear function in $\Co(\Delta)$ with $T_i(v_i)=1$
and $T_i(v_j)=0$ for $i\neq j$. 

\begin{Thm}[Billera \cite{MR1013666}, Thm.~2.3 and Thm.~3.6]\label{thm:Billera} The tent
  functions $T_1,\dots,T_m$ generate $\Co(\Delta)$. The kernel of the
  map $\R[X_1,\dots,X_m]\to\Co(\Delta)$ given by $X_i\mapsto T_i$ is
  generated by the elements
  \begin{itemize}
  \item $1-\sum_{i=1}^m X_i$
  \item $X_{i_1}\cdots X_{i_e}$ whenever
    $\{i_1,\dots,i_e\}\subset\{1,\dots,m\}$ is such that
    $\Delta(v_{i_1},\dots,v_{i_e})$ is not contained in $\Delta$.
  \end{itemize}
\end{Thm}
  
\begin{Remark}
  Billera makes the additional assumption that $\Delta$ be
  pure-dimensional. However, this appears to be immaterial
  for that particular part of his paper.
\end{Remark}

\smallskip We will mostly be concerned with the case when $\Delta$ is
of dimension $1$. We will write
$E_\Delta=\{(i,j)\in\{1,\dots,n\}^2\:|\: i<j\text{ and }
\Delta(v_i,v_j)\in\Delta\}$ for the set of indices corresponding to
the edges of $\Delta$. By Billera's theorem, $\Co(\Delta)$ is
generated by the tent functions $T_1,\dots,T_m$ subject to the rules:
\begin{enumerate}
\item $\sum_{i=1}^m T_i=1$.
\item $T_iT_jT_k=0$ for all distinct $i$,$j$,$k$.
\item $T_iT_j=0$ if and only if $(i,j)\notin E_\Delta$.
\end{enumerate}

We will have to go back and forth between the two descriptions of
$\Co(\Delta)$ that we have seen. Our description of the affine
variety associated with $\Delta$ in Prop.~\ref{prop:DirectLimit}, together with explicit
coordinates, translates into the following proposition, the proof of which
is obvious:

\begin{Prop}\label{prop:LimitLines}
  Let $\Delta\subset\R^n$ be a purely $1$-dimensional simplicial
  complex with vertices $V_\Delta=\{v_1,\dots,v_m\}$. For every edge
  $(i,j)\in E_\Delta$, let $C_{ij}$ be a copy of $\R$ with coordinate
  ring $\R[t_{ij}]$. Let $\phi$ be the unique map from the disjoint
  union of all $C_{ij}$, $(i,j)\in E_\Delta$, into $\R^n$ taking
  $C_{ij}$ to the line $v_i+\R\cdot (v_i-v_j)$ and mapping $-1\in
  C_{ij}$ to $v_i$ and $1\in C_{ij}$ to $v_j$. The dual ring
  homomorphism $\phi^\ast\colon\Co(\Delta)\to\prod_{(i,j)\in E_\Delta}
  \R[t_{ij}]$ induces an isomorphism
\[
\Co(\Delta)\isom\left\{(f_{ij})_{(i,j)\in E_\Delta}\:\biggl|\: 
  \begin{array}{ll}
f_{ij}(-1)=f_{ik}(-1) & \text{if } (i,j),(i,k)\in E_\Delta\\
f_{ij}(1)=f_{jk}(-1) & \text{if } (i,j),(j,k)\in E_\Delta\\    
f_{ij}(1)=f_{kj}(1) & \text{if }(i,j),(k,j)\in E_\Delta
  \end{array}; i,j,k=1,\dots,m
\right\}
\]
Under this isomorphism, the tent function $T_k$ corresponds to the
function $(f_{ij})_{(i,j)\in E_\Delta}$ with $f_{ik}=\frac 12(1+t_{ik})$ and
$f_{kl}=\frac 12(1-t_{kl})$ for all $(i,k),(k,l)\in E_\Delta$ and $f_{ij}=0$ for all
$(i,j)\in E_\Delta$ with $i,j\neq k$.

\noindent In particular, the product $T_iT_j$ with $(i,j)\in E_\Delta$ is
$(0,\dots,\frac 14(1-t_{ij}^2),\dots,0)$.\qed
\end{Prop}

We have to say what the degree of a piecewise polynomial should be:

\begin{Def} Let $\Delta=\sigma_1\cup\cdots\cup\sigma_k\subset\R^n$.
  Given $F\in\Co(\Delta)$, let $\sR(F)$ be the set of all $k$-tuples
  of polynomials $f_i\in\R[t_1,\dots,t_n]$ with
  $F|_{\sigma_i}=f_i|_{\sigma_i}$. Define the \emph{degree} of $F$,
  denoted $\deg(F)$, as
  $\deg(F)=\min_{(f_i)\in\sR(F)}\{\max_i\{\deg(f_i)\}$. 
\end{Def}

\begin{Remark} If we identify $\Co(\Delta)$ with
  $\bigl\{(g_1,\dots,g_k)\in\prod_{i=1}^k R[t_1,\dots,t_{d_i}]\
  \bigl|\ g_i|_{U_{ij}}=g_j|_{U_{ij}}\}$ as in
  Prop.~\ref{prop:DirectLimit}, then $F$ has a unique representation
  $F=(g_1,\dots,g_k)$ with $g_i\in\R[t_1,\dots,t_{d_i}]$, and
  $\deg(F)=\max_i\{\deg(g_i)\}$.

  On the other hand, every $F\in\Co(\Delta)$ can be expressed
  (non-uniquely) as a polynomial in the tent functions
  $T_1,\dots,T_m$. Let $\sT(F)=\{G\in\R[t_1,\dots,t_m]\;|\; F =
  G(T_1,\dots,T_m)\}$ and define $\deg_\sT(F)=\min\{\deg(G)\; |\;
  G\in\sT(F)\}$. Since $\deg(T_i)=1$, we have
  $\deg(F)\le\deg_\sT(F)$. In general, this inequality may be
  strict. For example, if $\Delta$ consists of two isolated points,
  then $\deg(F)=0$ for all $F\in\Co(\Delta)$ but $\deg_\sT(F)=1$
  whenever $F$ is non-constant.
\end{Remark}

\begin{Remark}\label{rem:LinearExtension}
  Write $\Co(\Delta)$ as in Prop.~\ref{prop:LimitLines} and fix an
  edge $(k,l)\in E_\Delta$. Then given any $g\in\R[t_{kl}]$, there
  exists $(f_{ij})\in\Co(\Delta)$ with $f_{kl}=g$ and $\deg(f_{ij})\le
  1$ for all $(i,j)\neq (k,l)$. In fact, we can take $f_{ij}=0$ for
  $\{i,j\}\cap\{k,l\}=\emptyset$, and if $(k,i)\in E$ for some $i\neq
  l$, put $f_{ki}=\frac{g(-1)}{2}(1-t_{ki})$, and similarly in the
  remaining cases.
\end{Remark}

Finally, we set up some notation and terminology for
quadratic modules: Let $A$ be a ring (commutative with unit). By $\sum
A^2$, we denote the set of all sums of squares of elements in $A$. For
finitely many elements $h_1,\dots,h_r\in A$, we write
\[
\QM(h_1,\dots,h_r)=\biggl\{s_0+\sum_{i=1}^r s_i h_i\ \bigl|\ s_0,\dots,s_r\in\sum
A^2\biggr\}
\]
and call this the \emph{quadratic module} generated by
$h_1,\dots,h_r$. The quadratic module $M=\QM(h_1,\dots,h_r)$ is called
\emph{archimedean} if for every $f\in A$ there exists $N\in\Z_+$ such
that $N+f\in M$. 

\begin{Prop}\cite[Cor.~5.2.4]{MR2383959}\label{Prop:ArchQM} With $A$ and $M$ as above,
  assume that $A$ is finitely
  generated over a field by elements $t_1,\dots,t_m$. The following
  are equivalent:
  \begin{enumerate}
  \item $M$ is archimedean.
  \item There exists $N\in\Z_+$ such that $N-\sum_{i=1}^m t_i^2\in M$.
  \item There exists $N\in\Z_+$ such that $N\pm t_i\in M$ for all $i=1,\dots,m$.
  \end{enumerate}
\end{Prop}

\section{Positivity in $\Co(\Delta)$}\label{sec:main}

Recall from the introduction the statement of Putinar's
positivstellensatz:

\begin{Thm}[Putinar \cite{MR1254128}]\label{thm:PutPPF} 
 Let $U$ be an affine $\R$-variety with coordinate ring $\R[U]$,
  let $h_1,\dots,h_r\in\R[U]$ and $K=\{x\in U\ |\ h_1(x)\ge
  0,\dots,h_r(x)\ge 0\}$. If the quadratic module $\QM(h_1,\dots,h_r)$
  is archimedean, then it contains every $f\in\R[U]$ such that $f(x)>0$
  holds for all $x\in K$.
\end{Thm}

In the original paper, as well as in \cite{MR1829790}, the theorem is
stated for $U=\R^n$. But it is straightforward to pass to the version
given here: Fix a surjection $\phi\colon\R[t_1,\dots,t_N]\to\R[U]$ and
a finite set of generators $G_1,\dots,G_s$ of $\ker(\phi)$, giving an
embedding $U=\{x\in\R^N\ |\ G_1(x)=\cdots=G_s(x)=0\}$. Choose
$H_1,\dots,H_r\in\R[t_1,\dots,t_N]$ such that $\phi(H_i)=h_i$. Then
$K=\{x\in\R^N\ |\ H_i(x)\ge 0, G_j(x)\ge 0, -G_j(x)\ge 0\text{ for all
}i,j\}$. Put $M=\QM(h_1,\dots,h_r)\subset\R[U]$ and
$M_0=\QM(H_1,\dots,H_r,\pm G_1,\dots,\pm G_s)\subset
R[t_1,\dots,t_n]$, so that $\phi(M_0)=M$. That $M$ is archimedean
means that there exists $N\in\Z_+$ such that $N-\sum \phi(t_i)^2\in
M$. From this we conclude $N-\sum t_i^2\in M_0$, so that $M_0$ is
archimedean, too. Given $f\in\R[U]$ as in the theorem, we may choose
$F\in\R[t_1,\dots,t_N]$ with $\phi(F)=f$ and conclude $F\in
M_0$. Applying $\phi$ gives the desired representation of $f$ in
$M$. Alternatively, one can deduce the above version of Putinar's
result directly from Jacobi's more general representation theorem (see
\cite{MR1838311} or \cite[Thm.~5.4.4]{MR2383959}).

\begin{Cor}
  Let $\Delta\subset\R^n$ be a simplicial complex with $m$
  vertices. Let $T_1,\dots,T_m\in\Co(\Delta)$ be the tent functions on
  $\Delta$. If a function $F\in\Co(\Delta)$ is strictly positive on
  $\Delta$, then there exist sums of squares
  $S_i$ in $\Co(\Delta)$, $i=0,\dots,m$, such that
\[
F= S_0 + \sum_{i=1}^m S_iT_i
\]
\end{Cor}

\begin{proof}
  Let $U$ be the affine variety associated with $\Delta$, as defined
  in Prop.~\ref{prop:DirectLimit}. Then it is easy to check that the
  tent functions define $\Delta$ as a subset of $U$,
  i.e.~$\Delta=\{x\in U(\R)\ |\ T_1(x)\ge 0,\dots,T_m(x)\ge 0\}$. The
  quadratic module $\QM(T_1,\dots,T_m)$ is archimedean by
  Prop.~\ref{Prop:ArchQM}, since it contains $1-T_i=\sum_{j\neq i}
  T_j$ for all $i=1,\dots,m$.
\end{proof}

\begin{Remark}
  Degree bounds for the sums of squares $S_i$ that depend only on the degree of $F$ cannot exist as soon as the dimension of $\Delta$ is at least two (see Scheiderer \cite{MR2182447}). However, there exist bounds that depend on other data, in particular the minimum of $F$ on $\Delta$ (see Schweighofer \cite{MR2068157}).
\end{Remark}

The following is a special case of the general results of Kuhlmann,
Marshall, and Schwartz for subsets of the line (see \cite{MR2174483}, \S4). 

\begin{Thm}\label{thm:KMSLine}
Every $f\in\R[t]$ such that $f|_{[-1,1]}\ge 0$ admits a representation
\[
f = s_0 + s_1(1-t^2)
\]
where $s_0,s_1$ are sums of two squares with $\deg(s_0)\le\deg(f)+1$, $\deg(s_1)\le\deg(f)-1$.
\end{Thm}

\begin{proof}
  Let $f\in\R[t]$ with $f|_{[-1,1]}\ge 0$. By
  \cite[Thm.~4.1]{MR2174483}, there is a representation
  $f=r_0+r_1(1+t)+r_2(1-t)+r_3(1-t^2)$ with the degree of each summand
  bounded by $\deg(f)$. Now substitute the identity $(1\pm t)=\frac
  12(1\pm t)^2+\frac 12(1-t^2)$.
\end{proof}

\noindent Translated into our setup, this says:

\begin{Cor}\label{cor:KMSLine}
  If $\Delta$ is the $1$-simplex with tent functions $T_1,T_2$, then
a function $F\in\Co(\Delta)$ is non-negative on $\Delta$ if and only
  if there exist sums of two squares $S, S_{12}$ in $\Co(\Delta)$ of
  degree at most $\deg(F)+1$ such
  that
\[
F=S + S_{12}T_1T_2
\]\qed
\end{Cor}

\noindent Our main result is a generalization to $1$-dimensional
simplicial complexes, which we restate from the introduction. 

\begin{Thm}\label{thm:Main}
  Let $\Delta\subset\R^n$ be a simplicial complex of dimension $1$
  with $e$ edges and $m$ vertices, of which $m_0$ are isolated. Let
  $T_1,\dots,T_m\in\Co(\Delta)$ be the tent functions on $\Delta$. A
  function $F\in\Co(\Delta)$ is non-negative on $\Delta$ if and only
  if there exist sums of squares $S, S_{ij}$ in $\Co(\Delta)$ such
  that
\[
F=S + \sum_{(i,j)\in E} S_{ij} T_iT_j.
\]
More precisely, there exist such $S$ and $S_{ij}$ of degree at
most $\deg(F)+6(e-1)+1$, where $S$ is a sum of at most $2e+m_0$
squares and each $S_{ij}$ is a sum of two squares.

\end{Thm}

\begin{Remark}
  The quadratic module $\QM(T_iT_j\,|\,(i,j)\in E)$ used in the theorem
  coincides in fact with the quadratic module $\QM(T_1,\dots,T_m)$
  used earlier. This follows from the identities $T_i=T_i(\sum_{j=1}^m
  T_j)=T_i^2+\sum_{(i,j)\in E} T_iT_j$ and $T_iT_j=T_iT_j(\sum_{k=1}^m
  T_k)=T_iT_j(T_i+T_j)=T_i^2T_j+T_j^2 T_i$, which was pointed out to
  me by the referee. In particular, $\QM(T_1,\dots,T_m)$ is in fact a
  preordering, i.e.~it is closed under multiplication. Using these
  identities, one could also restate the degree bounds in
  Thm.~\ref{thm:Main} for $\QM(T_1,\dots,T_m)$.
\end{Remark}

I would like to thank Claus Scheiderer for suggesting the
proof of the following lemma, replacing a much more pedestrian
argument in an earlier version:

\begin{Lemma}\label{lem:Gefusel}
  Let $f\in\R[t]$ be such that $f(x)\ge 0$ for all $x\in [-1,1]$.  For
  every $a,b\in\R$ with $a^2\le f(-1)$ and $b^2\le f(1)$, there exists
  $s\in\R[t]$ with $\deg(s^2)\le\deg(f)+3$ such that $s(-1)=a$,
  $s(1)=b$ and such that $s^2(x)\le f(x)$ for all $x\in [-1,-1]$.
\end{Lemma}

\begin{proof}
  By Thm.~\ref{thm:KMSLine}, there exist sums of squares
  $s_0,s_1\in\R[t]$ such that $f=s_0+s_1(1-t^2)$ and with $s_0$ of
  degree $2d\le\deg(f)+1$. Factor $s_0=\prod_{i=1}^d
  (t-\lambda_i)(t-\ol{\lambda_i})$ over $\C$ and let
  $g=\prod_{i=1}^d (t-\lambda_i)$. It follows that $g(-1)=\sqrt{f(-1)}\alpha$,
  $g(1)=\sqrt{f(1)}\beta$, where $\alpha,\beta\in\C$ with
  $|\alpha|=|\beta|=1$. Put
    \[
    \ell(t)=\frac{\ol\alpha a(1-t)}{2\sqrt{f(-1)}} + \frac{\ol\beta b(1+t)}{2\sqrt{f(1)}}
    \]
    and let $s=\re(g\cdot \ell)\in\R[t]$ (where the real part is taken on
    coefficients). By the choice of $\ell(t)$, we find that
    $s(-1)=a$ and $s(1)=b$. Furthermore, since $|\ell(x)|\le 1$ for
    all $x\in [-1,1]$, we obtain
\[
s^2(x)=\bigl(\re(g(x)\cdot\ell(x))\bigr)^2\le |g(x)\cdot\ell(x)|^2\le
|g(x)|^2=s_0(x)\le f(x)
\]
for all $x\in [-1,1]$.
\end{proof}

\begin{Cor}\label{cor:AdaptSOS}
  Let $f\in\R[t]$ be such that $f(x)\ge 0$ for all $x\in
  [-1,1]$. Given $k\in\Z_+$ and vectors $a,b,\in\R^k$ such that
  $f(-1)=||a||$ and $f(1)=||b||$, there exist polynomials
  $s_1,\dots,s_{k+2}\in\R[t]$ and $r\in\sum\R[t]^2$ such that the
  following hold:
  \begin{enumerate}
  \item $f=\sum_{i=1}^{k+2} s_i^2 + r(1-t^2)$.
  \item $s_i(-1)=a_i$, $s_i(1)=b_i$ for all $i=1,\dots,k$
  \item $s_i(-1)=0$, $s_i(1)=0$ for $i=k+1,k+2$
  \item $\deg(r),\deg(s_i^2)\le\deg(f)+3k+1$ for all $1\le i\le k+2$.
  \end{enumerate}
\end{Cor}

\begin{proof}
  Let $d=\deg(f)$. By the lemma, we may choose $s_1$ with
  $\deg(s_1^2)\le d+3$ such that $s_1(-1)=a_1$, $s_1(1)=b_1$ and
  $s_1^2\le f$, hence $f-s_1^2\ge 0$ on $[-1,1]$. Continuing
  inductively, we find $s_1,\dots,s_k$ such that $f-\sum_{i=1}^k
  s_i^2\ge 0$ holds on $[-1,1]$ and $s_i(-1)=a_i$, $s_i(1)=b_i$ for
  all $i=1,\dots,k$, with $\deg(s_i^2)\le d+3i$. By Thm.~\ref{thm:KMSLine}, there
  exist $s_{k+1},s_{k+2}$ and $r\in\sum\R[t]^2$ such that
\[
f-\sum_{i=1}^k s_i^2 = s_{k+1}^2 + s_{k+2}^2 + r\cdot (1-t^2)
\]
and with $\deg(s_{k+1}^2),\deg(s_{k+2}^2),\deg(r)\le d+3k+1$, which is
the desired representation. Note that condition (3) follows
automatically from (2), since $s_{k+1}^2(\pm 1)+s_{k+2}^2(\pm 1)=f(\pm
1)-\sum_{i=1}^k s_i^2(\pm 1) = 0$. 
\end{proof}

\begin{proof}[Proof of Thm.~\ref{thm:Main}]
  We do induction on the number $e$ of edges in $\Delta$. If $e=0$,
  then $\Delta$ is just the set $\{v_1,\dots,v_{m_0}\}$ of isolated
  vertices. In this case, $T_i=T_i^2$ for all $i$ and
  $F=\sum_{i=1}^{m_0} F(v_i)T_i$ is a sum of squares. 

  If $\Delta$ is not connected, say $\Delta=\Delta_1\cup\Delta_2$ with
  $\Delta_1\cap\Delta_2=\emptyset$ and
  $\Delta_1,\Delta_2\neq\emptyset$, then we can write
  $\Co(\Delta)=\Co(\Delta_1)\times\Co(\Delta_2)$. Applying the
  induction hypothesis to $\Delta_1$ and $\Delta_2$ gives the result.

  Now assume that $\Delta$ is connected. If $e=1$, the statement
  reduces to that of Corollary \ref{cor:KMSLine}, so we assume $e\ge
  2$ and $(1,2)\in E_\Delta$. Let $\Delta_1=\Delta(v_1,v_2)$,
  and let $\Delta_2$ be the closure of $\Delta\setminus\Delta_1$ so
  that $\Delta=\Delta_1\cup\Delta_2$ and
  $\Delta_1\cap\Delta_2\subset\{v_1,v_2\}$. We treat the case
  $\Delta_1\cap\Delta_2=\{v_1,v_2\}$. (If $\Delta_1\cap\Delta_2$
  contains only one vertex, the argument is analogous but somewhat simpler.)

Let $F\in\Co(\Delta)$ be non-negative on $\Delta$. By
Prop.~\ref{prop:LimitLines}, we can write $F=(f(t),F_2)$ with
$f\in\R[t]$ a polynomial in one variable satisfying $f(-1)=F_2(v_1)$
and $f(1)=F_2(v_2)$. By the induction hypothesis applied to $F_2$,
there is a representation
\[
F_2=\sum_{i=1}^{2(e-1)} \wt S_i^2 + \sum_{(i,j)\in E_{\Delta_2}} \wt R_{ij} T_i T_j.
\]
with $\wt R_{ij}$ sums of two squares in $\Co(\Delta_2)$, and with
$\deg(\wt S_i^2),\deg(\wt R_{ij})\le\deg(F)+6(e-2)+1$. By Cor.~\ref{cor:AdaptSOS}, we can write
\[
f=\sum_{i=1}^{2e} s_i^2 + r_{12}(1-t^2)
\]
such that $s_i(-1)=S_i(v_1)$ and $s_i(1)=S_i(v_2)$ for all
$i=1,\dots,2(e-1)$ and $s_i(-1)=s_i(1)=0$ for 
$i=2e-1,2e$, where $\deg(s_i^2)\le\deg(f)+6(e-1)+1$. It follows that
\[
S_i=\left\{
  \begin{array}{ll}
    (s_i,\wt S_i) & \text{ for } i=1,\dots,2(e-1)\\
    (s_i,0) & \text{ for }i=2e-1,2e
  \end{array}\right.
\]
are well-defined elements of $C^0(\Delta)$. Choose sums of two squares
$R_{ij}\in\sum\Co(\Delta)$, $(i,j)\in E$, such that
$R_{12}|_{\Delta_1}=r_{12}$, $R_{ij}|_{\Delta_2}=\wt R_{ij}$, and
$\deg(R_{12})=\deg(r_{12})$, $\deg(R_{ij})=\deg(\wt R_{ij})$, for all
$(i,j)\in E_{\Delta_2}$ (see Remark \ref{rem:LinearExtension}). Since $T_{ij}$ is
supported on $\Delta(v_i,v_j)$, we see that
\[
F=\sum_{i=1}^{2e} S_i^2 + \sum_{(i,j)\in E_\Delta} R_{ij}T_iT_j.
\]
\end{proof}

\begin{Examples}
  \begin{enumerate}
  \item Since the tent functions are themselves non-negative on
    $\Delta$, they must have a representation as in Theorem
    \ref{thm:Main}. This is reflected in the simple identity
    $T_i=T_i(\sum_{j=1}^m T_j)=\sum_{j=1}^m T_iT_j$.

  \item Let $\Delta$ be the boundary of the triangle spanned by
    $v_1=(0,1)$, $v_2=(0,0)$, $v_3=(1,0)$ in $\R^2$. The
    Zariski-closure of $\Delta$ is the plane curve $C=\{(x,y)\in\R^2\
    |\ xy(1-x-y)=0\}$, a union of three lines. We can write
    $\Co(\Delta)$ in terms of the tent functions $T_1,T_2,T_3$ with
    the relations $T_1+T_2+T_3=1$ and $T_1T_2T_3=0$, or we can
    identify it with the coordinate ring $\R[C]$ which is isomorphic
    to
\[
\bigl\{(f,g,h)\in \R[u]\times\R[v]\times\R[w]\ |\ f(1)=g(-1), g(1)=h(-1), h(1)=f(-1)\bigr\}
\]
in such a way that $2T_1=(u+1,1-v,0)$, $2T_2=(1-u,0,1+w)$,
$2T_3=(0,v+1,1-w)$ (see Prop.~\ref{prop:LimitLines}). Consider the
function $F=(u^2,v^2,w^2)=1-4T_1T_2-4T_1T_3-4T_2T_3$. From the second
expression, it is not immediately clear that $F$ is non-negative on
$\Delta$, while this is obvious from the first, since $F$ is even
non-negative on all of $C$. But $F$ is not a square in $\R[C]$, since
$(u,v,w)\notin\R[C]$, nor is it even a sum of squares in $\R[C]$ (see
\cite{MR2652535}, Example (1) in the introduction). However, by
Thm.~\ref{thm:Main}, it is contained in the quadratic module
\[
\QM(T_1T_2,T_1T_3,T_2T_3)=\QM((1-u^2,0,0),(0,1-v^2,0),(0,0,1-w^2)).
\]
Using the idea of the proof of Thm.~\ref{thm:Main}, one quickly arrives at the representation $F=(u^2,-v,w)^2+(u,-v,-1)^2(1-u^2,0,0)$. Translated into tent functions, this corresponds to the equality
\[
1-4T_1T_2-4T_1T_3-4T_2T_3 = (4T_1T_2-2T_1-2T_2+1)^2+(4T_1-2)^2T_1T_2.
\]
 \end{enumerate}
\end{Examples}

{
\small\linespread{1}

}

\begin{thebibliography}{10}
\providecommand{\url}[1]{\texttt{#1}}
\providecommand{\urlprefix}{URL }

\bibitem{MR1013666}
L.~J. Billera.
\newblock The algebra of continuous piecewise polynomials.
\newblock \emph{Adv. Math.}, \textbf{76~(2)}, 170--183, 1989.
\newline\urlprefix\url{http://dx.doi.org/10.1016/0001-8708(89)90047-9}

\bibitem{MR1838311}
T.~Jacobi.
\newblock A representation theorem for certain partially ordered commutative
  rings.
\newblock \emph{Math. Z.}, \textbf{237~(2)}, 259--273, 2001.
\newline\urlprefix\url{http://dx.doi.org/10.1007/PL00004868}

\bibitem{MR2174483}
S.~Kuhlmann, M.~Marshall, and N.~Schwartz.
\newblock Positivity, sums of squares and the multi-dimensional moment problem.
  {II}.
\newblock \emph{Adv. Geom.}, \textbf{5~(4)}, 583--606, 2005.

\bibitem{MR2383959}
M.~Marshall.
\newblock \emph{Positive polynomials and sums of squares}, vol. 146 of
  \emph{Mathematical Surveys and Monographs}.
\newblock American Mathematical Society, Providence, RI, 2008.

\bibitem{MR2110098}
E.~Miller and B.~Sturmfels.
\newblock \emph{Combinatorial commutative algebra}, vol. 227 of \emph{Graduate
  Texts in Mathematics}.
\newblock Springer-Verlag, New York, 2005.

\bibitem{MR2652535}
D.~Plaumann.
\newblock Sums of squares on reducible real curves.
\newblock \emph{Math. Z.}, \textbf{265~(4)}, 777--797, 2010.
\newline\urlprefix\url{http://dx.doi.org/10.1007/s00209-009-0541-8}

\bibitem{MR1829790}
A.~Prestel and C.~N. Delzell.
\newblock \emph{Positive polynomials}.
\newblock Springer Monographs in Mathematics. Springer-Verlag, Berlin, 2001.
\newblock From Hilbert's 17th problem to real algebra.

\bibitem{MR1254128}
M.~Putinar.
\newblock Positive polynomials on compact semi-algebraic sets.
\newblock \emph{Indiana Univ. Math. J.}, \textbf{42~(3)}, 969--984, 1993.

\bibitem{MR2020709}
C.~Scheiderer.
\newblock Sums of squares on real algebraic curves.
\newblock \emph{Math. Z.}, \textbf{245~(4)}, 725--760, 2003.

\bibitem{MR2182447}
---{}---{}---.
\newblock Non-existence of degree bounds for weighted sums of squares
  representations.
\newblock \emph{J. Complexity}, \textbf{21~(6)}, 823--844, 2005.

\bibitem{MR1092173}
K.~Schm{\"u}dgen.
\newblock The {$K$}-moment problem for compact semi-algebraic sets.
\newblock \emph{Math. Ann.}, \textbf{289~(2)}, 203--206, 1991.

\bibitem{MR2068157}
M.~Schweighofer.
\newblock On the complexity of {S}chm\"udgen's positivstellensatz.
\newblock \emph{J. Complexity}, \textbf{20~(4)}, 529--543, 2004.

\bibitem{MR0327060}
G.~Strang.
\newblock Piecewise polynomials and the finite element method.
\newblock \emph{Bull. Amer. Math. Soc.}, \textbf{79}, 1128--1137, 1973.

\end{thebibliography}
\end{document}